\documentclass[11pt]{article}
\usepackage{amsmath}
\usepackage{amsfonts}
\usepackage{amssymb}
\usepackage{enumerate}
\usepackage{amsthm}
\usepackage{graphicx}
\usepackage{verbatim}
\newtheorem{myproposition}{Proposition}[section]
\newtheorem{mytheorem}[myproposition]{Theorem}
\newtheorem{mylemma}[myproposition]{Lemma}

\newtheorem{mydefinition}[myproposition]{Definition}

\newtheorem{myobservation}[myproposition]{Observation}

\newtheorem{myproblem}[myproposition]{Problem}

\def\ni{\noindent}

\def\gr{\mathcal{G}}
\def\gA{\mathcal{A}}
\def\zet{\mathbb{Z}}

\makeatletter
\def\imod#1{\allowbreak\mkern10mu({\operator@font mod}\,\,#1)}
\makeatother

\makeatletter
\def\jmod#1{\allowbreak\mkern10mu{\operator@font mod}\,\,#1}
\makeatother

\begin {document}

\title{Note on group distance magic graphs $G[C_4]$}

\author{Sylwia Cichacz\thanks{The author was partially supported by National Science Centre grant nr 2011/01/D/ST/04104, as well as by the Polish Ministry of Science and Higher Education.}\\
{\small Faculty of Applied Mathematics}\\
{\small AGH University of Science and Technology}\\
{\small Al. Mickiewicza 30, 30-059 Krak\'ow, Poland}}

\maketitle

\begin{abstract}
A \emph{group distance magic labeling} or a $\gr$-distance magic labeling
of a graph $G(V,E)$ with $|V | = n$ is an injection $f$ from $V$ to an Abelian group $\gr$ of
order $n$ such that the weight $w(x)=\sum_{y\in N_G(x)}f(y)$ of every vertex $x \in V$ is equal to the same element
$\mu \in \gr$, called the magic constant.
In this paper we will show that if $G$ is a graph of order $n=2^{p}(2k+1)$ for some natural numbers $p$, $k$ such that  $\deg(v)\equiv c \imod {2^{p+1}}$ for some constant $c$ for any $v\in V(G)$, then there exists an $\gr$-distance magic labeling
for any abelian group $\gr$ for the graph $G[C_4]$. Moreover we prove that if $\gr$ is an arbitrary abelian group of order $4n$ such that $\gr \cong \zet_2 \times\zet_2 \times \gA$ for some abelian group $\gA$ of order $n$, then exists a $\gr$-distance magic labeling for any graph $G[C_4]$.\\

\end{abstract}


\section{Introduction}

All graphs considered in this paper are simple finite graphs. Consider a simple graph $G$ whose order we denote by $|G|=n$. Write $V (G)$ for the vertex set and $E(G)$ for the edge set of a graph $G$. The \emph{neighborhood} $N(x)$ of a vertex $x$ is the set of vertices adjacent to $x$, and the degree $\deg(x)$ of $x$ is $|N(x)|$, the size of the neighborhood of $x$. \\

Let $w(x)=\sum_{y\in N_G(x)}l(y)$ for every $x \in V(G)$.\\

\emph{Distance magic labeling} (also called \emph{sigma
labeling}) of a graph $G=(V,E)$ of order $n$ is a bijection $l
\colon V \rightarrow \{1, 2,\ldots  , n\}$ with the property that
there is a positive integer $k$ such
that
$w(x) = k$ for every $x \in V$. If a graph $G$ admits a distance magic labeling, then we say that $G$ is \emph{distance magic graph} (\cite{Vi}). The concept of distance magic labeling has been motivated by the construction of magic squares.\\

The following observations were independently proved:

\begin{myobservation}[\cite{Ji,MRS,Rao,Vi}] Let $G$ be a $r$-regular distance magic graph on $n$
vertices. Then $k = \frac{r(n+1)}{2}$.\end{myobservation}

\begin{myobservation}[\cite{Ji,MRS,Rao,Vi}] No $r$-regular graph with $r$-odd can be a distance
magic graph.\label{nieparzyste}\end{myobservation}

Problem of distance magic labeling of $r$-regular graphs was studied recently (see \cite{AFK,C,FKK1,MRS,RSP}).
It is interesting that if you blow up a $r$-regular $G$ graph into some specific $p$-regular graph (like $C_4$ or $\overline{K}_{2n}$), then the
obtained graph $H$ is distance magic. More formally, we have the following definition.
\begin{mydefinition}
Let $G$ and $H$ be two graphs where $\{x^1, x^2,\ldots,x^p\}$ are vertices
of $G$. Based upon the graph $G$, an isomorphic copy $H^j$ of $H$ replaces every vertex
$x^j$ , for $j = 1, 2,\ldots , p$ in such a way that a vertex in $H^j$ is adjacent to a vertex in
$H^i$ if and only if $x^jx^i$ was an edge in $G$. Let $G[H]$ denote the resulting graph.
\end{mydefinition}
Miller at al. \cite{MRS} proved the following results.
\begin{mytheorem}[\cite{MRS}] The cycle $C_n$ of length $n$ is a distance magic graph if and only if $n = 4$.
\end{mytheorem}

\begin{mytheorem}[\cite{MRS}]  If $r \geq 1$, $n \geq
3$, $G$ is an $r$-regular graph and $C_n$ the cycle of length $n$.
Then $G[C_n]$ admits a distance magic labeling if and only if $n =
4$.
\end{mytheorem}

\begin{mytheorem}[\cite{MRS}]  Let $G$ be an arbitrary regular graph. Then $G[\overline{K}_n]$ is distance
magic for any even $n$.
\end{mytheorem}

The following problem was posted  in \cite{AFK}.

\begin{myproblem}[\cite{AFK}] If $G$ is non-regular graph, determine if there is a distance magic labeling of
      $G[C_4]$.
\end{myproblem}

It seems to be very hard to characterize such graphs. For example there were considered all graphs $K_{m,n}[C_4]$  for $1\leq m<n\leq2700$ and only  $K_{9,21}[C_4]$, $K_{20,32}[C_4]$, $K_{428,548}[C_4]$ are distance magic (see  \cite{AC}).\\

Froncek in \cite{Fro} defined the notion of \emph{group distance magic graphs}, i.e. the graphs allowing the bijective labeling of vertices with elements of an abelian group resulting in constant sums of neighbor labels.

\begin{mydefinition}
A group distance magic labeling or a $\gr$-distance magic labeling
of a graph $G(V,E)$ with $|V | = n$ is an injection $f$ from $V$ to an abelian group $\gr$ of
order $n$ such that the weight $w(x)=\sum_{y\in N_G(x)}f(y)$ of every vertex $x \in V$ is equal to the same element
$\mu \in \gr$, called the magic constant.
\end{mydefinition}

 Obviously, every graph with $n$ vertices and a
distance magic labeling also admits a $\mathbb{Z}_n$-distance magic labeling. The converse is
not necessarily true.\\

It was proved that Observation~\ref{nieparzyste} is also true for $\gr$-distance magic labeling (\cite{CF}) in case $n\equiv 2(\mod 4)$.

\begin{myobservation}[\cite{CF}] Let $G$ be a $r$-regular distance magic graph on $n\equiv 2(\mod 4)$
vertices, where $r$ is odd.
There does not exists an abelian group $\gr$ of
order $n$ such that $G$ is $\gr$-distance magic.
\end{myobservation}

In this paper we will prove that if $G$ is a graph of order $n=2^{p}(2k+1)$ for some natural numbers $p$, $k$ such that  $\deg(v)\equiv c \imod {2^{p+1}}$ for some constant $c$ for any $v\in V(G)$, then there exists an $\gr$-distance magic labeling
for any abelian group $\gr$ for the graph $G[C_4]$.  Moreover we show that if $\gr$ is an abelian group of order $4n$ such that $\gr \cong \zet_2 \times\zet_2 \times \gA$ for some abelian group $\gA$ of order $n$, then there exists a $\gr$-distance magic labeling for any graph $G[C_4]$.

\section{Main results}
We start with the following lemma.

\begin{mylemma}\label{inne}
Let $G$ be a graph of order $n$ and $\gr$ be an arbitrary abelian group of order $4n$ such that $\gr \cong \zet_{2^{p}}  \times \gA$ for $p\geq 2$ and some abelian group $\gA$ of order $\frac{n}{2^{p-2}}$. If $\deg(v)\equiv c \jmod {2^{p-1}}$ for some constant $c$ and any $v\in V(G)$, then there exists a $\gr$-distance magic labeling for the graph $G[C_4]$.
\end{mylemma}

\ni\textbf{Proof.}
Let $G$ has the vertex set $V(G)=\{x^0, x^1,\ldots,x^{n-1}\}$, $C_4=v_0v_1v_2v_3v_0$ and $H=G[C_4]$.

For $0\leq i \leq n-1$ and $j=0,1,2,3$, let $v_j^i$ be the vertices of $H$ that replace
$x^i$, $0 \leq  i\leq n-1$ in $G$.

If $g \in \gr$, then we can write that $g=(w,a_i)$ for $w \in \mathbb{Z}_{2^p}$ and $a_i \in \gA$ for $i=0,1,\ldots,n-1$.

 Label the vertices of $H$ in the following way

$$f(v_j^i)=\left\{\begin{array}{lcc}
             \left((2i+j)\jmod {2^{p-1}},a_{\lfloor i \cdot 2^{-p+2}\rfloor}\right) & \mathrm{for} & j=0,1 \\
             \left(2^p-1,0\right)-f(v_{j-2}^i) & \mathrm{for} & j=2,3 \\
          \end{array}\right.$$
for $i=0,1,\ldots,n-1$. \\

 Notice that for every $i$
$$f(v_0^i) +f(v_2^i)=f(v_1^i) + f(v_3^i)=(2^{p}-1,0).$$
So the sum of the labels in the $i$th part is
$$f(v_0^i) + f(v_1^i)+f(v_2^i)+f(v_3^i) =(2^{p}-2,0),$$
which is independent of $i$. Since $c\equiv \deg(v) \imod {2^{p-1}}$ for any $v\in V(G) $, therefore, for every $x\in V(H)$

$$w(x)=(-2c-1,0).$$ \qed

\begin{mytheorem}\label{lemma22}
Let $G$ be a graph of order $n$ and $\gr$ be an arbitrary abelian group of order $4n$ such that $\gr \cong \zet_2 \times\zet_2 \times \gA$ for some abelian group $\gA$ of order $n$. There exists a $\gr$-distance magic labeling for the graph $G[C_4]$.
\end{mytheorem}
\ni\textbf{Proof.}
Let $G$ has the vertex set $V(G)=\{x^0, x^1,\ldots,x^{n-1}\}$  and $C_4=v_0v_1v_2v_3v_0$.
For $0\leq i \leq n-1$ and $j=0,1,2,3$, let $v_j^i$ be the vertices of $H$ that replace
$x^i$, $0 \leq  i\leq n-1$ in $G$.

If $g \in \gr$, then we can write that $g=(j_1,j_2,a_i)$ for $j_1,j_2 \in \mathbb{Z}_2$ and $a_i \in \gA$ for $i=0,1,\ldots,n-1$.

 Label the vertices of $H$ in the following way

$$f(v_j^i)=\left\{\begin{array}{ccc}
            (0,0,a_i) & \mathrm{for} & j=0, \\
            (1,0,a_i) & \mathrm{for} & j=1, \\
            (1,1,-a_i) & \mathrm{for} & j=2, \\
            (0,1,-a_i) & \mathrm{for} & j=3 \\
          \end{array}\right.$$
for $i=0,1,\ldots,n-1$.\\

Notice that for every $i=0,\ldots,n-1$
$$f(v_0^i) +f(v_2^i)=f(v_1^i) + f(v_3^i)=(1,1,0).$$
So the sum of the labels in the $i$th part is
$$f(v_0^i) + f(v_1^i)+f(v_2^i)+f(v_3^i) =(0,0,0),$$
which is independent of $i$. Therefore, for every $x\in V(H)$,
$$w(x) = (1,1,0).$$
\qed

\begin{mytheorem}\label{main}
Let $G$ be a graph of order $n$ and $\gr$  be an abelian group of order $4n$. If $n=2^{p}(2k+1)$ for some natural numbers $p$, $k$ and   $\deg(v)\equiv c \imod {2^{p+1}}$ for some constant $c$ for any $v\in V(G)$, then there exists a $\gr$-distance magic labeling for the graph $G[C_4]$.
\end{mytheorem}
\ni\textbf{Proof.}\\
The fundamental theorem of finite abelian groups states that the finite abelian group $\gr$ can be expressed as the direct sum of cyclic subgroups of prime-power order. This implies that $\gr \cong \mathbb{Z}_{2^{\alpha_0}} \times \mathbb{Z}_{p_{1}^{\alpha_1}}\times \mathbb{Z}_{p_{2}^{\alpha_2}}\times \ldots \times \mathbb{Z}_{p_{m}^{\alpha_m}}$ for some $\alpha_0>0$,
where $4n=2^{\alpha_0}\prod_{i=1}^m{p_i^{\alpha_i}}$ and $p_i$ for $i=1,\dots,m$ are not necessarily distinct primes. \\

\ni Suppose first that $\gr \cong \zet_2 \times \zet_2 \times \gA$ for some abelian group $\gA$ of order $n$, then we are done by Theorem~\ref{lemma22}.
Observe now that the assumption  $\deg(v)\equiv c \jmod {2^{p+1}}$  and unique decomposition of any natural number $c$ into powers of $2$ apply that  there exist constants $c_1,c_2,\ldots,c_{p}$  such that
$\deg(v)\equiv c_i \jmod {2^{i}}$  for $i=1,2,\ldots,p$ for any $v\in V(G)$. Hence if $\gr \cong \mathbb{Z}_{2^{\alpha_0}} \times \gA$ for some $2 \leq \alpha_0\leq p+2$ and some abelian group $\gA$ of order $\frac{4n}{2^{\alpha_0}}$, then
 we obtain by Lemma~\ref{inne} that there exists a $\gr$-distance magic labeling for the graph $G[C_4]$. \qed\\

The observation follows easily from the above Theorem~\ref{main}, however the below Observation~\ref{mlynek}  shows an infinite family of Eulerian graphs with odd order such that none of graphs was distance magic.
 \begin{myobservation}\label{drogi_zam}
Let $G$ be a graph of odd order $n$ and $\gr$  be an abelian group of order $4n$. If $G$ is an Eulerian graph (i.e. all vertices of the graph $G$ have even degrees), then there exists a $\gr$-distance magic labeling for the graph $G[C_4]$.
\end{myobservation}

Before we proof the  Observation~\ref{mlynek} we need the following definition. The \emph{Dutch windmill graph} $ D_m^{t}$ is the graph obtained by taking $t>1$ copies of the cycle $C_m$ with a vertex $c$ in common (\cite{Ga}). Thus for $t$ being even a graph $D_4^{t}$ is an Eulerian graph of odd order $3t+1$. Let
$i$th copy of a cycle in $C^{(t)}_4$ is $cy^ix^iz^ic$ for $i=0,\ldots,t-1$, $C_4=v_0v_1v_2v_3v_0$ and $H=C^{(t)}_4[C_4]$.
For $0\leq i \leq t-1$ and $j=0,1,2,3$, let $x^i_j$ ($y_j^i$, $z_j^i$, $c_j$ resp.) be the vertices of $H$ that replace
$x^i$ ($y^i$, $z^i$, $c$, resp.) $0 \leq  i\leq t-1$ in $C^{(t)}_3$.
\begin{myobservation}\label{mlynek}
There does not exist a distance magic graph $C^{(t)}_4[C_4]$.
\end{myobservation}
\ni\textbf{Proof.} Suppose that $C^{(t)}_4[C_4]$ is a distance magic graph. It is easy to observe that:
\begin{itemize}
  \item  $l(c_0)+l(c_{2})=l(c_1)+l(c_{3})=a_c$.
  \item  $l(x^i_0)+l(x^i_{2})= l(x^i_1)+l(x^i_{3}) = a_y^i$ for $0 \leq  i\leq t-1$.
  \item  $l(y^i_0)+l(y^i_{2})= l(y^i_1)+l(y^i_{3}) = a_y^i$ for $0 \leq  i\leq t-1$.
  \item $l(z^i_0)+l(z^i_{2})= l(z^i_1)+l(z^i_{3}) = a_z^i$ for $0 \leq  i\leq t-1$.
\end{itemize}

Since $w(z^i_j)=a_y^i+2a_x^i+2a_c=w(y^i_j)=a_z^i+2a_x^i+2a_c$ we obtain
that $a_y^i=a_z^i=a^i$ for $0\leq i \leq t-1$. Moreover $w(x^i_j)=a_x^i+4a^i=w(y^i_j)=a^i+2a_x^i+2a_c$, hence $a_x^i=3a_i-2a_c$.
Furthermore $w(z^i_j)=7a^i-2a_c=w(y^l_j)=7a^l-2a_c$ implies that $a^i=a^l=a$.

Since $7a-2a_c=a_c+4ta=k$ it has to be that $3a_c=(7-4t)a$ and therefore  $t=1$ (because $a_c,a>0$), a contradiction.~\qed\\

The following observation shows that inverse of the Theorem~\ref{main} is not true:

 \begin{myobservation}\label{dwudzielne}
Let $K_{p,q}$ be such complete bipartite  graph  that $p$ is even and $q$ is odd and $\gr$  be an abelian group of order $4(p+q)$. There exists a $\gr$-distance magic labeling for the graph $G[C_4]$.
\end{myobservation}
\ni\textbf{Proof.}\\
    If  $\gr \cong \zet_2 \times\zet_2 \times \gA$ for some abelian group $\gA$ of order $p+q$, then there exists a $\gr$-distance magic labeling for the graph $K_{p,q}[C_4]$ by Theorem~\ref{lemma22}.
    Suppose now that $\gr \cong \zet_4 \times \gA$ for some abelian group $\gA$ of order $p+q$.
Let $K_{p,q}$ has the partition vertex sets $A=\{x^0, x^1,\ldots,x^{p-1}\}$, $B=\{y^0, y^1,\ldots,y^{q-1}\}$  and $C_4=v_0v_1v_2v_3v_0$.
For $0\leq i \leq n-1$ and $j=0,1,2,3$, let $x_j^i$ ($y_j^l$ respectively) be the vertices of $K_{p,q}[C_4]$ that replace
$x^i$ $0 \leq  i\leq p-1$ ($y^i$ $0 \leq  l\leq q-1$ respectively) in $K_{p,q}$. If $g \in \gr$, then we can write that $g=(j,a_i)$ for $j \in \mathbb{Z}_4$ and $a_i \in \gA$ for $i=0,1,\ldots,p+q-1$.

 Label the vertices of $K_{p,q}[C_4]$ in the following way

$$f(x_j^i)=\left\{\begin{array}{ccc}
            (2j,a_i) & \mathrm{for} & j=0,1 \\
            \left(1,0\right)-f(x_{j-2}^i) & \mathrm{for} & j=2,3 \\
          \end{array}\right.$$
for $i=0,1,\ldots,p-1$.\\
$$f(y_j^l)=\left\{\begin{array}{ccc}
            (2j,a_{p+l}) & \mathrm{for} & j=0,1 \\
            \left(3,0\right)-f(y_{j-2}^l) & \mathrm{for} & j=2,3 \\
          \end{array}\right.$$
for $l=0,1,\ldots,q-1$.\\

Notice that
$$\begin{array}{c}
            f(x_0^i) +f(x_2^i)=f(x_1^i) + f(x_3^i)=(1,0),\\
            f(y_0^l) +f(y_2^l)=f(y_1^l) + f(y_3^l)=(3,0),\\
          \end{array}
$$

for every $i=0,\ldots,p-1$ and for every $l=0,\ldots,q-1$.
This implies:
$$\begin{array}{l}
\sum_{i=0}^{p-1}\left(\sum_{j=0}^3f(x_j^i)\right)=p(2,0)=(0,0)\\
\sum_{l=0}^{q-1}\left(\sum_{j=0}^3f(y_j^l)\right)=q(2,0)=(2,0)\end{array}$$
Hence  for every $x\in V(K_{p,q}[C_4])$,
$w(x)=(3,0)$. \qed\\


\end {document}